\documentclass[article,authoryear,jmlmc]{beg_32}             %  final version
\usepackage[hang]{footmisc}

\fancypagestyle{plain}{%
  \fancyhf{}
  \fancyhead[R]{\small {\it \jname}}
  \fancyfoot[R]{\small\bf\thepage }
  \fancyfoot[L]{\fottitle}
  }
\renewcommand{\dmyy}{20}

\begin{document}

\volume{}
\title{Transfer Learning on Multi-Fidelity Data}
\titlehead{Transfer Learning on Multi-Fidelity Data}
\authorhead{Dong H. Song \& Daniel M. Tartakovsky}
%For at least  authors with different addresses, use instead the following commands
\author[1]{Dong H. Song}
\corrauthor[1]{Daniel M. Tartakovsky}
% \author[2]{Third Author}
\corremail{tartakovsky@stanford.edu}
% \corraddress{Department of Energy Resources Engineering, Stanford University, Stanford, CA 94305, USA}
\address[1]{Department of Energy Resources Engineering, Stanford University, Stanford, CA 94305, USA}
% \address[2]{Business or Academic Affiliation 2, City, Province, Zip Code, Country Business or Academic Affiliation 2, City, Province, Zip Code, Country}
% End information for at least  authors with different addresses
% For authors with the same post address,
%\corrauthor{First A. Author}
%\corremail{f.author@affiliation.com}
%\author{Second B. Author, Jr.}
%\address{Department of Chemistry and Courant, Institute of Mathematical Sciences, New York, NY 10012, USA}
% End commands for all authors with the same address

\dataO{mm/dd/yyyy}
%\dataO{}
\dataF{mm/dd/yyyy}
%\dataF{}

\abstract{Neural networks (NNs) are often used as surrogates or emulators of partial differential equations (PDEs) that describe the dynamics of complex systems. A virtually negligible computational cost of such surrogates renders them an attractive tool for ensemble-based computation, which requires a large number of repeated PDE solves. Since the latter are also needed to generate sufficient data for NN training, the usefulness of NN-based surrogates hinges on the balance between the training cost and the computational gain stemming from their deployment. We rely on multi-fidelity simulations to reduce the cost of data generation for subsequent training of a deep convolutional NN (CNN) using transfer learning. High- and low-fidelity images are generated by solving PDEs on fine and coarse meshes, respectively. We use theoretical results for multilevel Monte Carlo to guide our choice of the numbers of images of each kind. We demonstrate the performance of this multi-fidelity training strategy on the problem of estimation of the distribution of a quantity of interest, whose dynamics is governed by a system of nonlinear PDEs (parabolic PDEs of multi-phase flow in heterogeneous porous media) with uncertain/random parameters. Our numerical experiments demonstrate that a mixture of a comparatively large number of low-fidelity data and smaller numbers of high- and low-fidelity data provides an optimal balance of computational speed-up and prediction accuracy. The former is reported relative to both CNN training on high-fidelity images only and Monte Carlo solution of the PDEs. The latter is expressed in terms of both the Wasserstein distance and the Kullback–Leibler divergence.
}

% OLD ABSTRACT

% Randomness in subsurface heterogeneity drives the uncertainty of subsurface flow characteristics. Current uncertainty quantification methods often involve repeated computationally expensive flow simulations. We train a deep convolutional neural network (CNN) surrogate model from multiple levels of data using transfer learning. The high fidelity data is of $[128 \times 128]$ and the intermediate fidelity data is of $[64 \times 64]$; the data is used to train the model in different phases, without modification. The model takes a high fidelity  permeability field as input and generates high fidelity saturation maps at time-steps of interest as output. The model is trained on a small budget (12 hours of simulation time) of multiple levels of data, and out performs models trained on 70+ hours of high fidelity data in terms of saturation map generation. The trained model is used to estimate the distributions (CDF/PDF) of water breakthrough times in multi-phase flow problems. In this uncertainty quantification task, the CNN surrogate performs comparable to 70+ hours of Monte Carlo simulations on various distribution metrics including Wasserstein Distance, and KL divergence.


\maketitle

 \graphicspath{
               {Figures/}
              }

\section{Introduction}
\label{sec:intro}
Machine learning techniques, especially neural networks (NNs), have pervaded every facet of human activity and has permeated into the field of scientific computing. %The NN based methods' ability to flexibly adapt to the training data and fast forward-pass run speeds are driving the network based method's popularity. 
In the latter setting, NNs are used to approximate highly nonlinear and irregular functions \citep{friedman2001elements}, solve (ordinary and partial) differential equations \citep[e.g.,][among many others]{lee1990neural, lagaris1998artificial, fuks2020}, and construct cheap surrogates for ensemble-based computation \citep[e.g.,][]{mo2019deep_b, raissi2019physics}. Examples of the latter include inverse modeling \citep{mo2019deep_a, zitong2020}, data assimilation \citep{tang2020deep}, and uncertainty quantification  \citep{tripathy2018deep, zhu2019physics}. 

A typical ensemble-based computation of practical significance involves repeated solves of (coupled, nonlinear) partial-differential equations (PDEs)
\begin{equation}
\label{eq:1}
    \mathcal N(\mathbf u; \boldsymbol{\theta}) = g(\mathbf x, t;\boldsymbol{\theta}), \qquad (\mathbf x,t) \in D \times (0,T],
\end{equation}
which describe the spatiotemporal evolution of (a set of) state variables $\mathbf u(\mathbf x,t)$ in the computational domain $D$ over simulation time horizon $(0,T]$. Multiple solves of~\eqref{eq:1}---for different values of the inputs $\boldsymbol{\theta}(\mathbf x,t)$ that parameterize the differential operator $\mathcal N$, the source function $g$, and auxiliary functions in the initial and/or boundary conditions---are required because these values are known at best in terms of their distributions, which are either inferred from data or provided by the expert. High computational cost of solving~\eqref{eq:1} numerically often precludes one from generating enough samples to obtain meaningful statistics of $\mathbf u(\mathbf x,t)$ or the derived quantities of interest. A surrogate of~\eqref{eq:1} carries a negligible cost, making possible ensemble-based computation with arbitrarily small sampling error. 

Alternative strategies for surrogate construction include polynomial chaos expansions \citep{xiu-2010-numerical}, Kriging or Gaussian processes \citep{Couckuyt-2014-ooDACE}, polynomial regression \citep{Montgomery-2018-Second}, tensor-product splines \citep{Hwang-2018-fast} and random forests \citep{Breiman-2014-Random}. Current popularity of NN-based surrogates \citep{mo2019deep_b, raissi2019physics} is grounded in the 
%
%Surrogate modeling using neural networks  enables the latter three applications by creating simple models. Such applications of surrogate modeling require many forward pass runs which can be computationally expensive on a traditional simulator. Surrogate models do not have to be NN based; data driven surrogate models can be made from other machine learning techniques including regression techniques and tree based techniques. However, NN models show promise based on their 
scalability and approximation capabilities of deep NNs \citep{friedman2001elements, tripathy2018deep}.
Regardless of the surrogate type, the training of a surrogate requires a large number of solves of~\eqref{eq:1} for different combinations of parameter values $\boldsymbol\theta$. Advanced computer architectures, e.g., CUDA-compatible graphics processing units (GPUs) and tensor processing units (TPUs), are almost a necessity to train a large NN. A combined cost of training-data acquisition and  NN training can be so large as to negate the benefits of the NN. 

This observation suggests that the practical utility of a NN as a surrogate model hinges on one's ability to dramatically reduce the cost of its construction. We rely on multi-fidelity simulations to reduce the cost of data generation for subsequent training of a deep convolutional NN (CNN) using transfer learning. High- and low-fidelity images are generated by solving PDE~\eqref{eq:1} on fine and coarse meshes, respectively.
%
%Traditionally, a simulation was built and run to solve the model numerically. A numerical simulation is a modeled system where the domain has been discretized (volumes, elements, etc) and the PDE's have been translated into a system of algebraic equations on the discretization. Real systems are often expressed by nonlinear models; a numerical linear solver must be employed to solve a nonlinear set of equations which is computationally expensive. 
A fine mesh is defined by the need to resolve the spatiotemporal variability of the model's inputs $\boldsymbol\theta$ and outputs $\mathbf u$; the resulting high-fidelity simulation carries a high computational cost. Lower-fidelity solutions of~\eqref{eq:1}, obtained on coarser meshes on which appropriately homogenized inputs $\boldsymbol\theta_\text{hom}$ are defined, are cheaper to compute but less accurate. %Furthermore, real systems are often heterogeneous; a good numerical simulation should be discretized to a degree of granularity which captures the effects of heterogeneity. A higher computation cost is associated with a finer mesh as there are more calculations are needed to solve the simulation. Also, the mesh density effects the accuracy of the numerical solution; a finer mesh will have less numerical error than a coarse mesh. The computational cost associated with repeated model forward pass can be prohibitive when traditional simulation is used for sample based applications such as UQ \citep{zitong2020}.
We train a CNN on a mixture of these multi-fidelity data, using the theoretical results for multilevel Monte Carlo (MLMC) \citep{heinrich1998monte, heinrich2001multilevel, giles2008multilevel, taverniers2020accelerated} to guide our choice of the numbers of solutions $\mathbf u(\mathbf x,t)$ of each kind.
%Surrogate models can greatly enhance the utility of sample based applications where the simulation costs are prohibitive. Compared to solving a PDE based model using a traditional simulator, a NN surrogate model is fast because the output is directly calculated without the need for a linear solver; a NN surrogate calculates $u(x,t)$ from Equation \ref{eq:1} using simple linear algebra multiplications and summations once the model weights have been trained. NN based surrogate models show strong potential to increase the utility of applications which traditionally required simulators to repeatedly solve PDE's.
Similar to MLMC \citep{muller-2013-Multifidelity, Peherstorfer-2019-Multifidelity}, the varying fidelity (aka ``levels'') of predictions of $\mathbf u$ can be achieved not only by solving~\eqref{eq:1} on different meshes, but also by replacing~\eqref{eq:1} with its cheaper-to-compute counterparts. For example, the multi-phase flow equations used as the computational testbed in this study can be replaced with the cheaper-to-solve Richards equation and Green-Ampt equation \citep{yang-2020-resource, sinsbeck-2015-impact}, each of which encapsulates progressively simplified physics. We leave this aspect of NN training on multi-fidelity data for a follow-up study.

Section~\ref{sec:Model} contains a brief description of our CNN and the workflow for its training  on multi-fidelity of data. The performance of this algorithm is tested on a system of nonlinear parabolic PDEs governing multi-phase flow in a heterogeneous porous medium with uncertain properties, which are formulated in Section~\ref{sec:ComputationalExample}. In Section~\ref{sec:Results}, we demonstrate  the accuracy and computational efficiency of the CNN-based surrogate used to quantify predictive uncertainty of~\eqref{eq:1} in terms of the distribution of a quantity of interest. Main conclusions drawn from this study are presented in Section~\ref{sec: Conclusions}.

%%%%%%%%%%%%%%%%%%%%%%%%%%%%%%%%%%%%%%%%%%%%%%%%%
\section{Deep Convolutional Neural Networks} 
\label{sec:Model}
%%%%%%%%%%%%%%%%%%%%%%%%%%%%%%%%%%%%%%%%%%%%%%%%%

While many flavors of NNs can be used as a surrogate for a PDE-based model like~\eqref{eq:1}, we choose CNNs because of their proven ability to model complex nonlinear phenomena and the negligible cost of their forward pass. To be concrete, we select the CNN with encoder-decoder architecture \citep{mo2019deep_a}, which has previously been used for single-phase \citep{mo2019deep_b} and multi-phase \citep{mo2019deep_a} flow problems in the context of uncertainty quantification. The encoder-decoder architecture is ideally suited for training on multi-fidelity data, as detailed in  Section \ref{subsec: transfer}.
%
%A deep CNN is selected as the platform for our surrogate model to meet the following requirements:
%\begin{itemize}
%    \item Ability to model complex-nonlinear processes
%    \item Very fast forward pass times
%    \item Able to be trained on multiple fidelity of data (HFS and LFS)
%\end{itemize}
%For this study, we did wanted to develop a method where we do not have to modify the training data of the surrogate model. A modified CNN surrogate model found in \cite{mo2019deep_a} is selected as the encoder-decoder architecture of the model met the above criteria. The original authors applied this network to solve single phase \cite{mo2019deep_b} and multi phase problems \cite{mo2019deep_a} for UQ; the fast forward pass speed of this model type is critical to the UQ analysis.   The network used in this study is described in Table~\ref{tab: model layers}:

The CNN-based surrogate is set up as an image-to-image regression model \citep{zitong2020}. To  train and test the network, we use the parameter values $\boldsymbol\theta(\mathbf x_i)$ in $N_\text{el}$ elements $\{ \mathbf x_i \}_{i = 1}^{N_\text{el}}$ of a numerical grid as  input and the discretized solution $\mathbf u(\mathbf x_i,t_k)$ of PDE~\eqref{eq:1} at $N_\text{ts}$ time steps $\{t_k\}_{k=1}^{N_\text{ts}}$ as output. To facilitate the generalizibility of the trained CNN to unseen sets of the input $\boldsymbol\theta(\mathbf x_i)$, i.e., to ensure that the CNN is not over-fitted to a particular choice of $\boldsymbol\theta(\mathbf x_i)$, the training data comprises a large number $N_\text{train}$ of the solutions $\mathbf u$ obtained for $N_\text{train}$ realizations $\{\boldsymbol \theta_1,\dots, \boldsymbol \theta_{N_\text{train}}\}$ of the input $\boldsymbol\theta$. The loss function, 
\begin{equation}\label{eq: loss function}
  \mathcal L (\mathbf w) = \sum_{m=1}^{N_\text{train}} \sum_{i=1}^{N_\text{el}} \sum_{k=1}^{N_\text{ts}} | \mathbf u(\mathbf x_i,t_k; \boldsymbol\theta_m) - \hat{\mathbf u}_{ik}(\mathbf w; \boldsymbol\theta_m)|+\lambda \sum_{n=1}^{N_\text{w}} w_n^2
\end{equation}
consists of two parts. The first represents the $L_1$-norm discrepancy between the state variables $\mathbf u$ predicted by solving PDE~\eqref{eq:1}, $\mathbf u(\mathbf x_i,t_k)$ and estimated by the CNN, $\hat{\mathbf u}_{ik}(\mathbf w)$, with $N_\text{w}$ weights $\mathbf w = (w_1,\dots,w_{N_\text{w}})^\top$. The $L_2$-norm  regularization term prevents over-fitting by penalizing large weights $\mathbf w$ associated complex models; the regularization parameter $\lambda$ determines how much regularization penalty is applied. The CNN training consists of finding a set of weights $\mathbf w^\star$ that minimizes $\mathcal L$.

\subsection{Transfer Learning}
\label{subsec: transfer}

The construction of CNN-based generalizable surrogates, which are capable of making predictions for realizations of $\boldsymbol\theta(\mathbf x)$ not seen during training, requires a large number of PDE solves, $N_\text{train}$; e.g., $N_\text{train} \sim 1500$ was used by \cite{mo2019deep_b} and \cite{zitong2020} to train the encoder-decoder CNNs similar to ours.  If a single PDE solve is expensive, the costs associate with large $N_\text{train}$ can be large to the point where CNN training becomes unfeasible. To alleviate this problem, we use both multi-fidelity data and transfer learning \citep{donahue2014decaf}. The latter is a technique that uses a NN trained for one task as the starting point for a different NN being trained for a new task.
%{\color{blue} The number of training data $N_\text{train}$ required to train an accurate and generalized CNN can be large especially if the output responses are strongly nonlinear \citep{mo2019deep_b}; to train a similar encoder-decoder CNN surrogate \cite{zitong2020} used 1600 simulation runs and \cite{mo2019deep_b} needed 1000-1500 simulation runs. In practice, it is common to use a NN pre-trained on readily available data as a starting point to minimize data and training costs. The final NN would be refined using data directly applicable to the new task. } 
Transfer learning has been implemented for face detection \citep{jiang2017face}, generation of image description \citep{karpathy2015deep}, and construction of physics-informed NNs \citep{haghighat2021physics}, among other applications. 

%\textbf{[rewrite in terms relevant to the above, what is transfer learning too, Add NLFS]}  Training large neural networks for complex tasks requires large amounts of data, and sufficient data may be costly or nonexistent. In practice, it is common to initially train a network on large and readily available data-sets similar to the data needed to train the network for the actual task. This trained network would be the starting point or a part of the final model which would be refined using fewer data directly applicable to the task. This concept is called transfer learning, and it was introduced by \citep{donahue2014decaf}. Transfer learning has been implemented for complex tasks such as face detection \citep{jiang2017face}, the generation of image description \citep{karpathy2015deep}, and the construction of physics-informed NNs \citep{haghighat2021physics}. 

\begin{figure}[htbp]
\begin{center}
\includegraphics[width=1\textwidth]{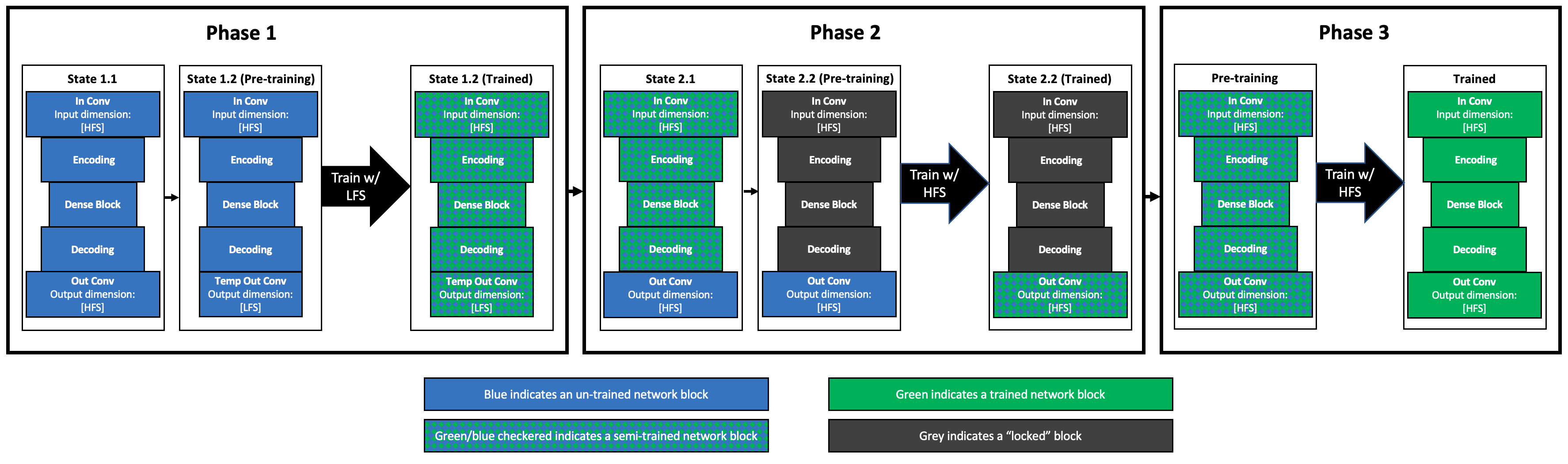}
\end{center}
\caption{Workflow for CNN training on multi-fidelity data. Phase 1 returns a low-resolution CNN trained on the LFS data. Phase 2 supplements that network with an additional layer whose weights are determined from the HFS data, producing a high-resolution CNN. In Phase 3, the latter is fine-tuned by allowing all the weights to vary during the training on the same HFS data. \ref{appendix: A} provides a pseudocode for all three Phases. }
\label{fig:method}
\end{figure}

Let HFS and LFS data refer to the solutions of~\eqref{eq:1}, $\mathbf u(\mathbf x_i, t_k)$, obtained on the fine ($N_\text{el} = N_\text{el}^\text{HFS}$) and coarse ($N_\text{el} = N_\text{el}^\text{LFS}$ with $N_\text{el}^\text{LFS} < N_\text{el}^\text{HFS}$) meshes, respectively. %Transfer learning allows a model trained using the LFS, the low cost readily available data, to be the starting point when training a model in the high-fidelity scale. 
If $N_\text{w}$ in~\eqref{eq: loss function} denotes the number of weights in the CNN trained on the HFS data, then our implementation of transfer learning starts with the construction of a CNN composed of $N_\text{LFS}$ ($N_\text{LFS}<N_\text{w}$) weights $\mathbf w_{\text{LFS}} = (w_1,\dots,w_{N_\text{LFS}})^\top$ trained on the LFS data. Then, the HFS data are used to train the desired high-resolution CNN, i.e., to determine the remaining weights $\mathbf w_{\text{HFS}} = (w_{N_\text{LFS}+1},\dots,w_{N_\text{w}})^\top$. This transfer learning strategy is depicted in Fig.~\ref{fig:method} and detailed below. %The resultant trained CNN, which performs on the fine scale, should have lower training data costs than a CNN trained using solely HFS data.

\subsection{Workflow for CNN Training on Multi-fidelity Data}
\label{sec:methodology}

Our strategy for CNN training on multi-fidelity data consists of three phases (Fig.~\ref{fig:method}), each of which results in a CNN denoted by $M_i$ ($i=1,2,3$). During Phase 1, the CNN $M_1$ with the $N_\text{el}^\text{LFS} \times N_\text{el}^\text{LFS}$ output is trained on the LFS data. In Phase 2, the CNN $M_2$ with $N_\text{el}^\text{HFS} \times N_\text{el}^\text{HFS}$ output is constructed by adding an additional layer with the weights $\mathbf w_\text{HFS}$, which are trained on the HFS data while keeping the original weights $\mathbf w_\text{LFS}$ fixed. Phase 3 consists of fine-tuning the CNN $M_2$ by allowing all the weights $\mathbf w = \{ \mathbf w_\text{LFS}, \mathbf w_\text{HFS} \}$ to update during the training on the same HFS data. The numerical experiments reported in Sections~\ref{sec:ComputationalExample} and~\ref{sec:Results} demonstrate that this transfer learning strategy significantly reduces the number of high-resolution PDE solves.

The workflow of our approach is provided below (see \ref{appendix: A} for the corresponding pseudocode).
%\textbf{why do we do this? as oposed to just training on coarse (after phase 1) how we handle we have a mismatch of low reso vs high reso}
\begin{enumerate}[leftmargin=2cm]
    \item[Phase 1:] Train a CNN $M_1$, with $N_\text{el}^\text{LFS} \times N_\text{el}^\text{LFS}$ output, on the LFS data.
    %{\color{blue} Train a CNN $M_1$ which is composed of $\mathbf w_{\text{LFS}} = (w_1,\dots,w_{N_\text{LFS}})^\top$ and a temporary layer $L_\text{temp}$. $M_1$ produces output of $N_\text{el}^\text{LFS} \times N_\text{el}^\text{LFS}$ and is trained on LFS data.}
    \begin{enumerate}
        \item[State 1.1:] Initialize the transfer learning by employing the encoder-decoder CNN of \cite{mo2019deep_a}, $M_\text{init}$ with  $N_\text{el}^\text{HFS} \times N_\text{el}^\text{HFS}$ output, whose $N_\text{w}$ weights $\mathbf{w}$ are set to \texttt{PyTorch} defaults.
        %{\color{blue} $M_\text{init}$ is disassembled according to $\mathbf w_{\text{LFS}}$ and $\mathbf w_{\text{HFS}}$. $\mathbf w_{\text{LFS}}$ is used as the starting point of $M_{ 1}$ described in State 2. %$\mathbf w_{\text{HFS}}$ is put aside for now as $L_{\text{last}}$.}
        \item [State 1.2:] Train the CNN $M_1$ on the LFS data. The starting point is a CNN obtained from $M_\text{init}$ by replacing its last layer $L_{\text{last}}$, which has $N_\text{HFS}$ weights $\mathbf w_{\text{HFS}}$, with a temporary convolution layer $L_{\text{temp}}$. The latter makes the output of $M_1$ match the dimensions of the LFS data, $[N_\text{ts} \times N_\text{el}^\text{LFS} \times N_\text{el}^\text{LFS}]$, and has significantly fewer weights than has  $L_{\text{last}}$. Then, $N_\text{LFS}$ weights $\mathbf w_{\text{LFS}}$ are updated by minimizing \eqref{eq: loss function} on the LFS data. 
    \end{enumerate}
    \item[Phase 2:] Train a CNN $M_2$, with $N_\text{w}$ weights $\mathbf w = \{ \mathbf w_\text{LFS}, \mathbf w_\text{HFS} \}$ (of which $N_\text{LFS}$ weights are locked) and $N_\text{el}^\text{HFS} \times N_\text{el}^\text{HFS}$ output, on the HFS data 
    %Train weights $\mathbf w_{\text{HFS}} = (w_{N_\text{LFS}+1},\dots,w_{N_\text{w}})^\top$ by minimizing \eqref{eq: loss function} using HFS data
    \begin{enumerate}
     %   \item[State 3:] $M_{1, best}$; low resolution outputs. $M_1$ trained using LFS
        \item[State 2.1:] Build a CNN from $M_1$ by replacing its layer  $L_{\text{temp}}$ with the layer $L_{\text{last}}$.
        \item[State 2.2:] Train the resulting CNN $M_2$ on the HFS data by minimizing~\eqref{eq: loss function} over the weights $\mathbf w_{\text{HFS}}$ of layer $L_{\text{last}}$, while keeping the remaining weights $\mathbf w_{\text{LFS}}$ fixed at their values in $M_1$.
    \end{enumerate}
    \item[Phase 3:] Train a CNN $M_3$ on the HFS data by allowing all weights $\mathbf w$ of $M_2$ to vary during the minimization
%    \begin{enumerate}
%        \item [State 6: ] $M_{2, best}$; high resolution outputs. $M_{2}$ trained using HFS.
 %       \item[State 7:] $M_{ 3}$; high resolution outputs. $M_{ 2, best}$ with all weights open for update during training.
  %      \item[State 8:] $M_{ 3, best}$; high resolution outputs. $M_{ 3}$ trained using HFS.
%    \end{enumerate}
\end{enumerate}
%{\color{blue} $M_{ 3, best}$ produces results in the high fidelity scale which was trained on LFS and HFS.} 
Since the bulk of the CNN $M_3$ training is carried out on the LFS data, this procedure is more efficient than CNN training solely on HFS data.  

%%%%%%%%%%%%%%%%%%%%%%%%%%%%%%%%%%%%%%%%%%%%%%%%%
\section{Computational Example: Multi-phase Flow}
\label{sec:ComputationalExample}
%%%%%%%%%%%%%%%%%%%%%%%%%%%%%%%%%%%%%%%%%%%%%%%%%

Numerical solution of problems involving multi-phase flow in porous media is notoriously difficult because of the high degree of nonlinearity and stiffness of the governing PDEs. Each forward solve of these PDEs is so expensive that it is uncommon, e.g., in petroleum engineering, to base uncertainty quantification efforts on as few as three model runs. %Applications which need repeated simulation runs, such as UQ and generating data for surrogate models, can have prohibitive computational costs. As such, the effectiveness of training a surrogate model on multiple fidelity  can be demonstrated through a this multi-phase flow computational example. 
This high cost and numerical complexity make the multi-phase flow equations a challenging testbed for ensemble-based simulations. 

We consider horizontal flow of two incompressible and immiscible fluids, with viscosities $\mu_1$ and $\mu_2$, in a heterogeneous, incompressible, and isotropic porous medium $D$. The latter is characterized by porosity $\phi$ and intrinsic permeability $k$. The porosity is assumed to be constant $\phi=0.25$, and intrinsic permeability $k(\mathbf x)$ is treated as a random variable. Mass conservation of the $\ell$th fluid phase ($\ell=1,2$) implies
\begin{subequations}\label{eq:2 mat bal}
\begin{equation}
    \phi \frac{\partial S_\ell}{\partial t} + \nabla\cdot \mathbf{v}_\ell + q_\ell = 0, \qquad \mathbf{x}\equiv (x_1,x_2)^\top \in D, \quad t \in [0,T],
\end{equation}
where $S_\ell(\mathbf{x},t)$ is the phase saturation constrained by $S_1+S_2=1$; %$D$ is the 2-dimensional spatial domain; $\phi$ is the porosity; and 
$q_\ell$ is the source/sink term; and the macroscopic velocity $\mathbf{v}_\ell(\mathbf{x},t)$ is described by the generalized Darcy law
\begin{align}\label{eq:3 darcy vel}
  \mathbf{v}_\ell = - k  \frac{k_{r\ell}}{\mu_\ell} \nabla P_\ell.
\end{align}
\end{subequations}
The relative permeability for the $\ell$th phase, $k_{r\ell}$, varies with the phase saturation, $k_{r\ell} = k_{r\ell}(S_\ell)$, in accordance with the Brooks-Corey constitutive model \citep{corey1954interrelation}. Following \cite{taverniers2020accelerated} and many others, we neglect the capillary forces, i.e., assume pressure within the two phases to be equal,  $P_1 = P_2 \equiv P(\mathbf{x},t)$; that is a common assumption in applications to reservoir engineering and carbon sequestration. %In this specific numerical example the subscripts $\ell$ = 1 and 2 represent water and oil respectively.

\begin{figure}[htbp]
\begin{center}
\includegraphics[width=0.5\textwidth]{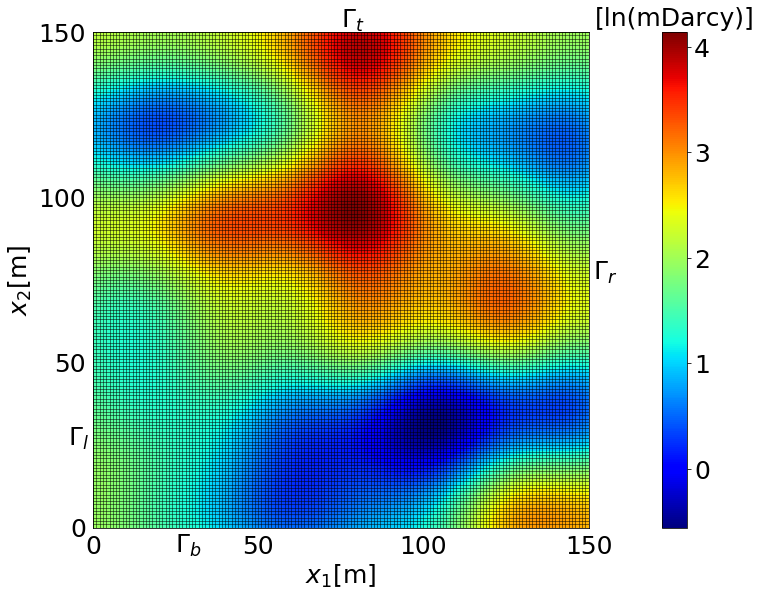}
\end{center}
\caption{A representative realization of log permeability field $Y=\ln k$ on the $128\times 128$ grid, which is used in high-fidelity simulations. Permeability $k$ is expressed in mDarcy.}
\label{fig: log_perm}
\end{figure}

The computational domain $D$ is a $150~\text{m} \times 150~\text{m}$ square (Fig.~\ref{fig: log_perm}) with the impermeable bottom ($\Gamma_\text{b}$ or $x_2 = 0$) and top ($\Gamma_\text{t}$ or $x_2 = 150$~m) boundaries; Dirichlet conditions are imposed along the left ($\Gamma_\text{l}$ or $x_1 = 0$) and right ($\Gamma_\text{r}$ or $x_1 = 150$~m) boundaries:
\begin{subequations}\label{eq:bc}
\begin{align}
    & \frac{\partial P}{\partial x_2} = 0, \quad \mathbf x \in \Gamma_\text{b} \cup \Gamma_\text{t}; \quad 
    P = 10.2 \;\;\&\;\; S_1 = 1.0, \quad \mathbf x \in \Gamma_\text{l}; \quad 
    P = 10.1, \quad \mathbf x \in \Gamma_\text{r}; 
 %   & S_1 = 1, \quad \mathbf x \in \Gamma_\text{l}
\end{align}
here and below, the pressure $P$ is expressed in MPa. Initial conditions are
\begin{align}
   P(\mathbf x, 0) = 10.1, \qquad S_1(\mathbf x, 0) = 0, \qquad \mathbf x \in D.
\end{align}
\end{subequations}
 %The top and bottom boundaries ( and  respectively) are no-flow boundaries (homogeneous Neumann conditions). 
%The left and right boundaries have Dirichlet conditions: the left boundary ($\Gamma_l$) is held at constant pressure $P=10.2$ MPa and saturation $S_1=1$, and the right boundary ($\Gamma_r$) is held at constant pressure $P=10.1$ MPa.

All the model parameters, except for the intrinsic permeability $k(\mathbf x)$, are assumed to be constant and known with certainty. The uncertain permeability $k(\mathbf x)$ is modeled as a second-order stationary random field, such that $Y(\mathbf{x}) =\ln k$ is multivariate Gaussian with mean $\langle Y \rangle = 0$, variance $\sigma^2_Y=2.0$, and an exponential two-point covariance $C(\mathbf{x},\mathbf{y})=\sigma^2_Y\exp (-|\mathbf{x}-\mathbf{y}|/\lambda_Y)$ with the correlation length $\lambda_Y=19$~m. %The coefficient of variation (CV) of $k=\exp(Y)$ is
%\begin{equation}\label{eq: coef var}
%\text{CV}(k) \equiv \frac{\sigma_k}{\langle k \rangle} = \frac{\sqrt{[\exp(\sigma_Y^2)-1]\exp(2\mu_Y + \sigma_Y^2)}}{\exp[\mu_Y + (\sigma_Y^2/2)]} = 2.53.
%\end{equation}
We use a truncated Karhunen-Lo\'eve expansion with $p=31$ terms to represent  $Y(\mathbf{x})$ \citep{taverniers2020accelerated}. A representative realization of the resulting permeability field is shown in Fig.~\ref{fig: log_perm} for the $128 \times 128$ mesh.

Equations~\eqref{eq:2 mat bal}--\eqref{eq:bc} are approximated using a finite volume scheme in space and implicit Euler scheme in time, yielding a highly nonlinear algebraic system \citep{aziz1979petroleum}.  The latter is solved, at each time step, through Newton-Raphson (NR) iterations with the modified Appleyard update dampening \citep{appleyard1981special} that improves the convergence of NR iterations by capping the maximum saturation update to a specified limit. %: in this computational example, $|S_{\ell,i}^{(\nu+1)} - S_{\ell,i}^{(\nu)}| \leq 0.3$ is used for the $\ell$ th phase, $\nu$ th NR iteration, at the $i$ th volume. 
For the $\nu$th iteration and the $i$th cell of volume $V_i$, the convergence criteria are
\begin{equation}\label{eq: conv crit}
  \max_i \left|\Delta t\left(\frac{r_{\ell,i}}{\phi V_i}\right)\right|<\epsilon_1, \quad
  \max_i |P_i^{(\nu+1)} - P_i^{(\nu)}|<\epsilon_2,
  \quad
  \max_i |S_{\ell,i}^{(\nu+1)} - S_{\ell,i}^{(\nu)}|<\epsilon_3
\end{equation}
where $r_{\ell,i}$ is the residual of the mass balance of phase $\ell$, $\Delta t$ is the time step, the relative residual norm $\epsilon_1 =10^{-6}$, the maximum pressure update $\epsilon_2 =10^{-3}$, and the maximum saturation update $\epsilon_3 =10^{-2}$.

\subsection{Upscaling of Permeability}

Multi-fidelity data are generated by solving~\eqref{eq:2 mat bal}--\eqref{eq:bc} on progressively coarsened grids: the $128 \times 122$ and $64 \times 64$ grids are used for HFS and LFS, respectively. %Surrogate models trained on multiple fidelity of data require simulation results at multiple fidelity. 
This grid coarsening must be accomplished by upscaling (coarsening) of the realizations of the random permeability $\hat{k}$ which are initially generated at the finest scale (Fig.~\ref{fig: log_perm}).  Among alternative upscaling strategies \citep{paleologos-1996-effective, tartakovsky-1998-transient-2, boso-2018-information}, we select the one proposed by \cite{durlofsky2005upscaling} because of its computational simplicity. It turns a scalar permeability field defined on the fine $(128 \times 128)$ mesh into its upscaled tensorial (anisotropic) counterpart whose off-diagonal components are 0 and the diagonal components are computed as the distance-weighted arithmetic mean perpendicular to the direction of flow and the distance-weighted harmonic mean in the direction of flow.
% Solving the computational example on the fine and the coarse scale permeability fields yields the training data at multiple fidelity. The original simulations were discretized into $128 \times 128$ grids to form the high fidelity simulations (HFS). The up-scaling scheme is used to generate the $64 \times 64$ intermediate fidelity simulations (LFS) and the $32 \times 32$ low fidelity simulations (LFS).

\subsection{Data Acquisition}

\begin{figure}[htbp]
\begin{center}
\includegraphics[width=0.8\textwidth]{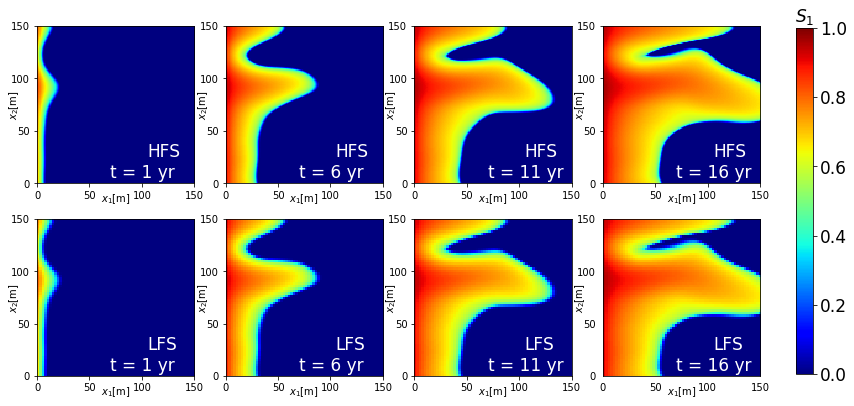}
\end{center}
\caption{Temporal snapshots of saturation $S_1(\mathbf x,t)$ computed with HFS  (top-row) and LFS (bottom-row) for the permeability field $k(\mathbf x)$ in Fig.~\ref{fig: log_perm}.}
\label{fig: saturations_data}
\end{figure}

Multi-fidelity training data come in the form of $N_\text{ts} = 16$ temporal snapshots of the saturation $S_1(\mathbf x,t)$ computed by solving~\eqref{eq:2 mat bal}--\eqref{eq:bc} on the $N_\text{el} \times N_\text{el}$ grids with $N_\text{el} = 128 \equiv N_\text{el}^\text{HFS}$ and $64 \equiv N_\text{el}^\text{LFS}$. Fig.~\ref{fig: saturations_data} shows examples of such images, corresponding to the permeability field in Fig.~\ref{fig: log_perm}. The permeability fields on the finest mesh, $[1 \times N_\text{HFS} \times N_\text{HFS}]$, are used as the input $\boldsymbol\theta$ for all CNNs. The size of of the CNN, $[N_\text{ts} \times N_\text{el} \times N_\text{el}]$, depends on the size of the training data. 
%Thus, each data set contains the input-output components:
%{\color{blue}
%\begin{itemize}
%    \item Permeability field [$1 \times N_\text{HFS} \times N_\text{HFS}$]
%    \item Saturation map at each time step [$N_\text{ts} \times N_\text{el} \times N_\text{el}$]
%\end{itemize}
%}

%{\color{blue} Comment this?} The upscaled permeability fields and pressure maps are not used to train the surrogate model. \textbf{[Why not?] {\color{blue} In this project we wanted to keep the input dimensions the same}} However, the former is used to check the integrity of the data sets. \textbf{[What does that mean?]}{\color{blue}Ask Daniel Should I just say the results of the simulations were used to train the models? also we probably dont need sentence}

%The results of the computational example are packed into data sets to train surrogate models. The specific computational example is a 2-dimensional simulation, so the data is treated like sets of [$nx_1 \times nx_2$] images as shown in Figs \ref{fig: log_perm} and \ref{fig: saturations_data} where a sample permeability field and associated saturation maps are displayed respectively. 
% Each data set contains three or four components depending on if the data set is an upscaled data set: \textbf{[I don't get it]}

% \begin{itemize}
%     \item Original permeability field  [$1 \times N_\text{HFS} \times N_\text{HFS}$]
%     \item Upscaled permeability field (only in upscaed data sets) [$1 \times nx_1 \times nx_2$]
%     \item Saturation map at each time step [$N_\text{ts} \times nx_1 \times nx_2$]
%     \item Pressure map at each time step [$N_\text{ts} \times nx_1 \times nx_2$]
% \end{itemize}
% where $[nx_{1i} \times nx_{2i}] = [128 \times 128]$ is the dimensions of the finest scale,  $[nx_1 \times nx_2] = [128 \times 128]$ (HFS) or $[64 \times 64]$ (LFS) is the dimensions of the computational example. 

The numerical solutions of~\eqref{eq:2 mat bal}--\eqref{eq:bc} are obtained using a \texttt{Matlab}-based multi-phase flow simulator on a computer with an Intel Core i7-4790 3.6GHz processor and 64GB of RAM. The computation time for each HFS data point is 219.13~sec and 37.13~sec for each LFS data point. The time needed to generate a data set is henceforth referred as ``data-generation budget''.

\subsection{CNN Training}

Table~\ref{tab: model layers} describes the CNN architecture used in this implementation our general approach (see Fig~\ref{fig:method}). The training is done with \texttt{PyTorch}, on the Stanford Mazama high-performance computing cluster. The allocated computing resources include Intel Xenon Gold 6126 CPU (2.6 GHz), 60GB RAM, and Nvidia V100 GPU with 16GB vRAM. (Although available, multi-cores were not used for this work.) %The associated hyper-parameter optimization is discussed in Subsection~\ref{subsec: hyperparameter}.

%The CNN is trained on the HFS and LFS data using transfer learning to achieve the best model performance for a given data budget. The transfer learning methodology shown in Fig~\ref{fig:method} include blocks which represent parts of the CNN; the actual layers and the associated input and output dimensions are provided in Table~\ref{tab: model layers}. LFS data are used to train $\mathbf w_{\text{LFS}}$ which compose the layers starting from Convolution 1 to Dense Block (Decoding). HFS data are used to train  $\mathbf w_{\text{HFS}}$ which compose the Convolution Transpose 2 layer. 

\begin{table}[htbp]
\caption{Model block description and the input and output dimensions of each model block. 
%In this table: $N_\text{ts}=16$, $N_\text{el}^\text{HFS} = 128$, $N_\text{el}^\text{LFS} = 64$, and $N_\text{dense} = 32$. The channels in each block are described by $n_c,c=[1,\dots,7]$ \textbf{[what are the values?] {\color{blue} [64,344,172,652,326,606,303]}}.
In our numerical experiments, the number of time steps is $N_\text{ts}=16$; the number of elements in fine and coarse meshes is $N_\text{el}^\text{HFS} = 128$ and $N_\text{el}^\text{LFS} = 64$, respectively; the number of elements in the output of the dense block is $N_\text{dense} = 32$; and the number of channels in each of the seven layers of the CNN is $n_1 = 64$, $n_2 = 344$, $n_3 = 172$, $n_4 = 652$, $n_5 = 326$, $n_6 = 606$, and $n_7 = 303$.}
\centering
\begin{tabular}{l c c}
\hline
 Layer  & Input & Output  \\
\hline
  Input: Permeability field $k$ & \multicolumn{2}{c}{$1 \times N_\text{HFS} \times N_\text{HFS}$}  \\
  Convolution 1  & $n_1 \times N_\text{el}^\text{HFS} \times N_\text{el}^\text{HFS}$ & $n_2 \times N_\text{el}^\text{LFS} \times N_\text{el}^\text{LFS}$  \\
  Dense Block (Encoding)  & $n_2 \times N_\text{el}^\text{LFS} \times N_\text{el}^\text{LFS}$ & $n_3 \times N_\text{el}^\text{LFS} \times N_\text{el}^\text{LFS}$   \\
  Convolution 2  & $n_3 \times N_\text{LFS} \times N_\text{el}^\text{LFS}$ & $n_4 \times N_\text{el}^\text{dense} \times N_\text{dense}$  \\
  Dense Block  & $n_4 \times N_\text{dense} \times N_\text{dense}$ & $n_5 \times N_\text{dense} \times N_\text{dense}$   \\
  Convolution Transpose 1  & $n_5 \times N_\text{dense} \times N_\text{dense}$ & $n_6 \times N_\text{el}^\text{LFS} \times N_\text{el}^\text{LFS}$    \\
  Dense Block (Decoding)  & $n_6 \times N_\text{el}^\text{LFS} \times N_\text{el}^\text{LFS}$ & $n_7 \times N_\text{el}^\text{LFS} \times N_\text{el}^\text{LFS}$  \\
  Convolution Transpose 2  & $n_7 \times N_\text{el}^\text{LFS} \times N_\text{el}^\text{LFS}$ & $N_\text{ts} \times N_\text{el}^\text{HFS} \times N_\text{el}^\text{HFS}$  \\
  Output: Saturation map $\hat S$ & \multicolumn{2}{c}{$N_\text{ts} \times N_\text{el}^\text{HFS} \times N_\text{el}^\text{HFS}$}  \\
\hline
% \multicolumn{2}{l}{$^{a}$Footnote text here.}
\end{tabular}
\label{tab: model layers}
\end{table}

The key hyper-parameters affecting the CNN performance are the learning rate (LR), the weight decay (WD), the factor (F), and the minimum learning rate (mLR). The LR and WD are parameters of the Adam  optimizer \citep{kingma2014adam}, and the F and mLR are parameters of the \texttt{ReduceLROnPlateau} scheduler. The CNN training involves many more hyper-parameters, but we use their default values in \texttt{PyTorch}. Further information on the hyper-parameters, schedulers, and optimizers can be found in the \texttt{PyTorch} documentation \citep{NEURIPS2019_9015}.

\begin{figure}[htbp]
\begin{center}
\includegraphics[width=1\textwidth]{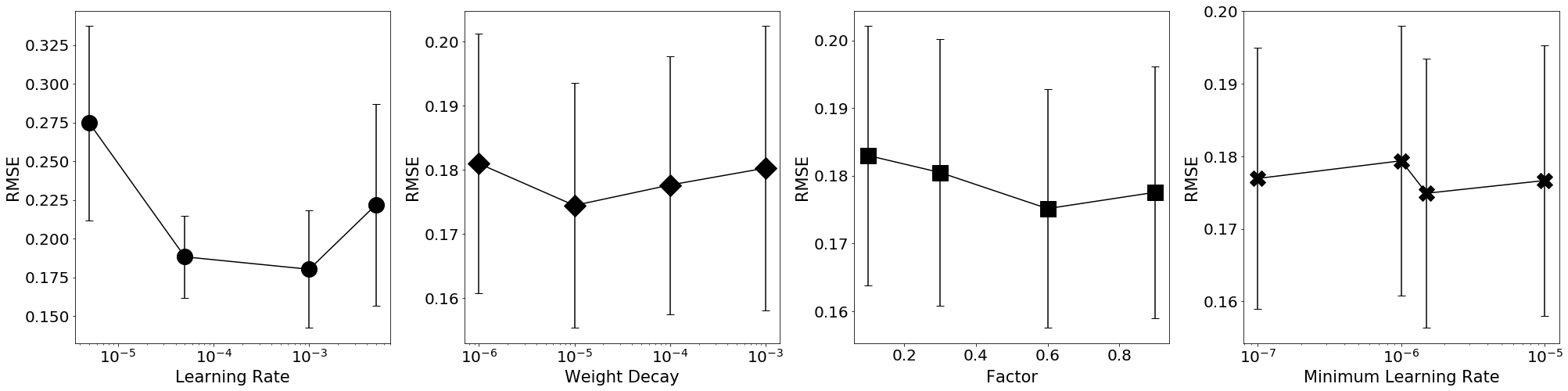}
\end{center}
\caption{Hyper-parameter performance in the neighborhood of optimum hyper parameter set in terms of the root mean square error (RMSE) for the test data. Unless labeled as the $x$-axis variable, all plot correspond to $\text{LR} = 5 \cdot 10^{-5}$, $\text{WD} =1 \cdot 10^{-5}$, $\text{F} =0.6$, and  $\text{mLR} = 5 \cdot 10^{-6}$. Each data point represents the mean and standard deviation of 10 training sessions.}
\label{fig: hyperparameter search}
\end{figure}

\begin{table}[htbp]
\caption{Learning rates and epochs used at each Phase.}
\centering
\begin{tabular}{cc c}
\hline
   & Learning rate & Epochs\\
\hline
  Phase 1  & $5 \cdot 10^{-4}$ & 170   \\
  Phase 2  & $5 \cdot 10^5$ & 150  \\
  Phase 3  & $10^{-5}$ & 100  \\
\end{tabular}
\label{tab:hyperparameters}
\end{table}

The hyper-parameters used by \cite{mo2019deep_a} in a similar CNN architecture serve as an initial guess for the hyper-parameter optimization. The latter required 100 HFS, with each training pass taking about 0.65 hours to complete, when 200 epochs were used. It took 7.2 training-hours to find optimal hyper-parameters (12 training passes), and a considerably smaller wall-clock time because this task was parallelized across several GPU nodes. We selected the hyper-parameter values yielding the smallest root mean square error (RMSE) on the HFS test data (Fig.~\ref{fig: hyperparameter search}). These values are used as a starting point in the hyper-parameter optimization for multi-fidelity transfer learning. Then, the LR and epochs at each Phase (Section~\ref{sec:methodology}) are modified to minimize the RMSE on the corresponding test data. The resulting hyper-parameter values are shown in Table~\ref{tab:hyperparameters}. 

\section{Results}
\label{sec:Results}

\begin{figure}[htbp]
\begin{center}
\includegraphics[width=1\textwidth]{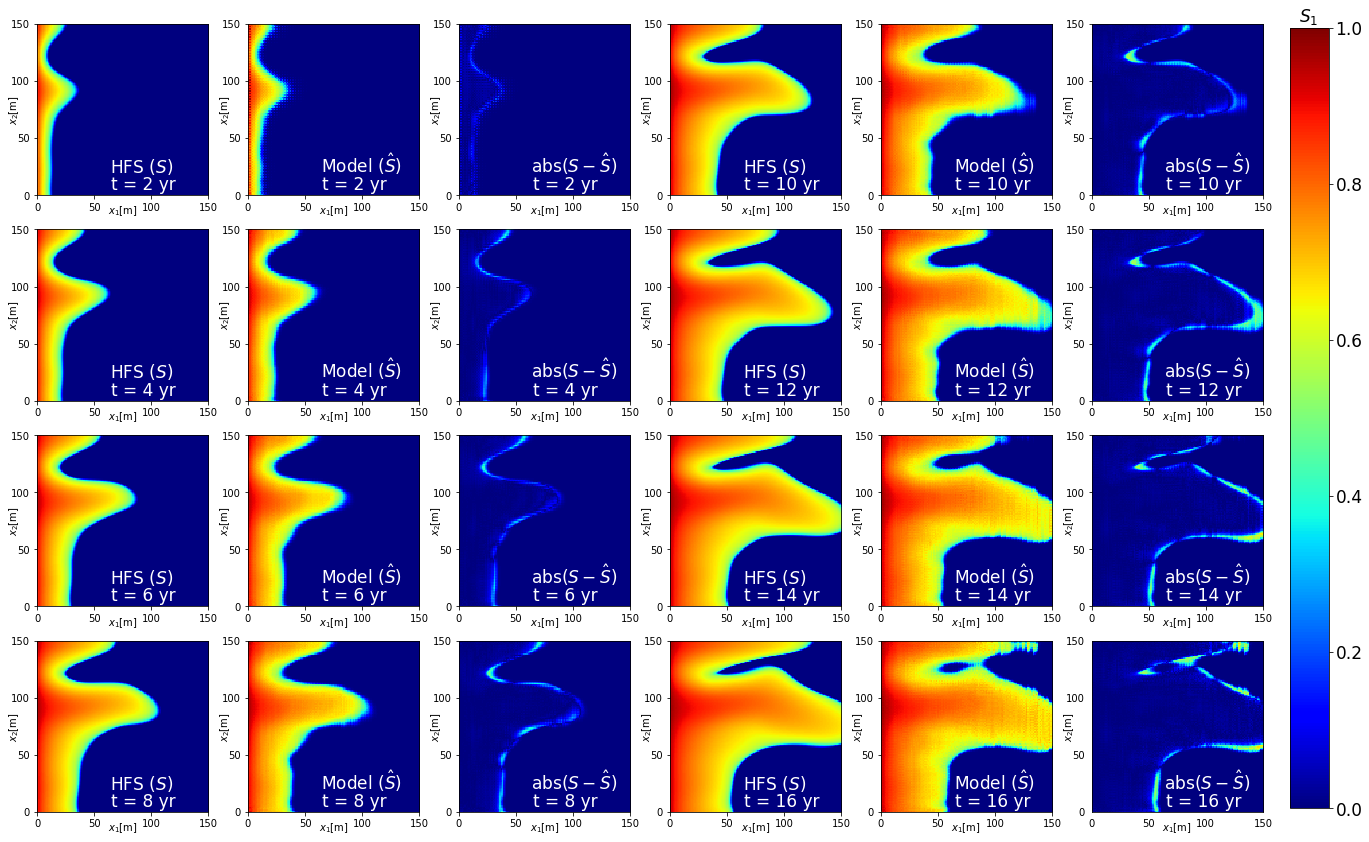}
\end{center}
\caption{Temporal snapshots of the saturation maps $S_1(\mathbf x,t)$ for the permeability field $k(\mathbf x)$ in Fig.~\ref{fig: log_perm}. These are generated with either HFS of the PDE model~\eqref{eq:2 mat bal} and~\eqref{eq:bc} (labeled as $S$ in the first and fourth columns) or the CNN surrogate (labeled as $\hat{S}$) in the second and fifth columns). The third and sixth columns display the absolute difference between the two predictions, $|S-\hat S |$.}
\label{fig: results diff}
\end{figure}

Once trained (in this example, on 573 LFS and 100 HFS, which took 12 hours to generate), the CNN surrogate provides an accurate approximation of the PDE solution on the fine mesh (Fig.~\ref{fig: results diff}), even for such highly nonlinear problems as~\eqref{eq:2 mat bal} that exhibit sharp dynamic fronts.  A forward pass of the CNN surrogate is on the order of a second, whereas a fine-mesh PDE solution takes nearly 220 seconds. This two-orders of magnitude speed up makes CNN surrogates an invaluable tool for UQ (Section~\ref{sec:UQ performance}). 

\subsection{Model Performance }
\label{sec:Model_performance vs ML Dat Ratio}

We compare the relative performance of the CNN trained on multi-fidelity data and the CNNs trained on either HFS data or LFS data, in terms of both accuracy (RMSE on test data) and computational cost. We also investigate the effect of varying the amount of HFS and LFS data for a given computational budget of 12 hours.

\begin{figure}[htbp]
\begin{center}
\includegraphics[width=\textwidth]{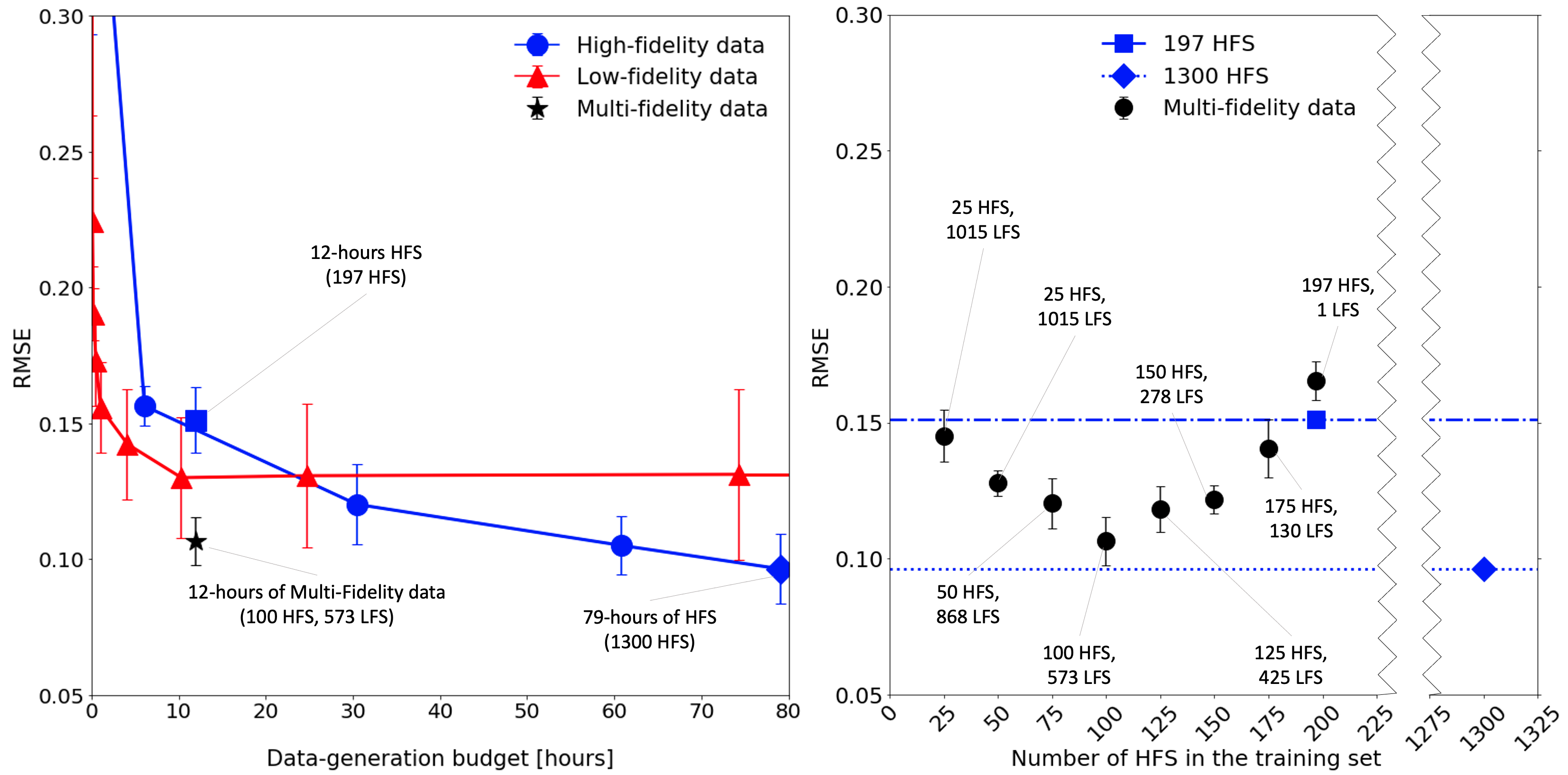}
\end{center}
\caption{RMSE on test data for the alternative CNN training strategies. It is plotted as function of the budget allocated for data generation (left) and the number of PDE solves on the fine mesh used to generate HFS data (right). Each RMSE point in these graphs represents an average over 10 iterations of training and is accompanied by error bars (the standard deviation). The left plate provides RMSE for the CNNs trained on high-fidelity (blue circles), low-fidelity (red triangles), or multi-fidelity (black star) data. The latter corresponds to the CNN trained on an optimal (the lowest RMSE) mix of high- and low-fidelity data for a set  budget of 12 hours; it is contrasted with the RMSE of the CNN trained on the HFS data generated within the same budget (blue square). The black circles in the right plate represent RMSE of the CNN trained on the multi-fidelity data sets, in which the number of HFS varies while the data-generation budget is fixed at 12 hours. Also shown there are RMSEs of the CNNs trained on 12 hours (dot-dashed line)  and 79 hours (dotted line) of HFS. }
\label{fig: results mix data}
\end{figure}

To train the high-resolution ($128 \times 128$ output) CNN solely on the LFS ($64 \times 64$) data, the latter have to be downscaled to match the dimensions. %simulation results had to be modified before training the network as the CNN surrogate model output is $[16 \times 128 \times 128]$ while the LFS output is $[16 \times 64 \times 64]$. This modification is achieved 
We do so by taking the Kronecker product of a $64 \times 64$ LFS image and a $2 \times 2$ matrix of 1s. The transformed LFS data have the desired dimensions, while containing the same information as the original image. The test data are composed of HFS images (PDE solves on fine mesh) that were not used for CNN training. %Training the CNN surrogate on the HFS is straightforward. 
Figure~\ref{fig: results mix data} exhibits the RMSEs on test data of the CNNs trained on high-, low-, and multi-fidelity data as function of the computational budget; each point in these graphs represents an average over 10 repetitions of training and is accompanied by error bars (the standard deviation).

The left plate of Fig.~\ref{fig: results mix data} reveals that, if the data-generation budget does not exceed 20~hours, the CNN trained on the LFS data outperforms its HFS-trained counterpart in terms of RMSE. That is because such budgets do not allow for generation of sufficient amounts of HFS data. As the budget increases, the error of the LFS data precludes RMSE of the CNN trained on such data from dropping below 0.125 while RMSE of the  HFS-trained CNN continues to decrease. This finding is reminiscent of the cost-constrained selection between high- and low-fidelity models in the context of ensemble-based simulations  \citep{yang-2020-resource, sinsbeck-2015-impact}. This figure also demonstrates that, for a relatively small budget of 12 hours, the use of multi-fidelity data yields the CNN whose RMSE is appreciably smaller that those of the CNNs trained on either HFS data or LFS data.

An optimal mix of the HFS and LFS data is investigated in the right plate of Fig.~\ref{fig: results mix data}. The multi-fidelity training was conducted 5 times for each HFS/LFS ratio, with random selection of LFS/HFS from a larger pool data. 
%A high data cost is associated with high performing surrogate models trained on HFS data. 
%When the CNN surrogates were trained using 100 HFS and 573 LFS, a data cost of 12 hours, the model is able to perform as models trained on 60-79hours of data cost.  
At the optimal mix of 573 LFS and 100 HFS, two of the five experiments, the CNN trained on 12 hours of these multi-fidelity data has lower RMSE than the mean RMSE of the CNN trained on 79 hours of the HFS data. This LFS/HFS ratio lies near the range, $1.5 - 5.5$, suggested for multilevel Monte Carlo \citep{taverniers2020accelerated}. For the data-generation budget of 12 hours, %The optimum training ratio is somewhere between almost all LFS and almost all HFS. When the data is almost all 
a mix dominated by the LFS data results in a CNN whose RMSE on test data exceeds 1.0, which indicates that the network's last Convolution Transpose 2 layer is not meaningfully trained. 

%Initially, the hyper parameter search for the CNN training on multi-level data was conducted at a data ratio nearest to the optimum point. A different set of hyper-parameters may change the optimum ratio. \textbf{[Chicken \& egg? {\color{blue} We can remove this}]}

%Although test RMSE is a common metric to determine the performance of models; a model with a low test RMSE does not necessarily mean the model has utility. The model is applied to a UQ task of calculating breakthrough times and the UQ results were compared to that of MC simulations. The UQ performance is discussed in Subsection~\ref{sec:Model_performance vs ML Dat Ratio}. 

\subsection{CNN Surrogates for Uncertainty Quantification}
\label{sec:UQ performance}

Finally, we investigate the utility of our CNN surrogates for uncertainty quantification. A quantity of interest is the breakthrough time, $T_\text{break}$, at the  $x_1 = 100$~m plane (Fig.~\ref{fig: log_perm}), with the term ``breakthrough'' defined as the saturation of the invading phase ($S_1$) exceeding 0.15. Given uncertainty in intrinsic permeability $k(\mathbf x)$, a solution of~\eqref{eq:2 mat bal} and, hence, predictions of $T_\text{break}$ are given in terms of their cumulative distribution functions (CDFs) or probability density functions (PDFs).

\begin{figure}[htbp]
\begin{center}
\includegraphics[width=1\textwidth]{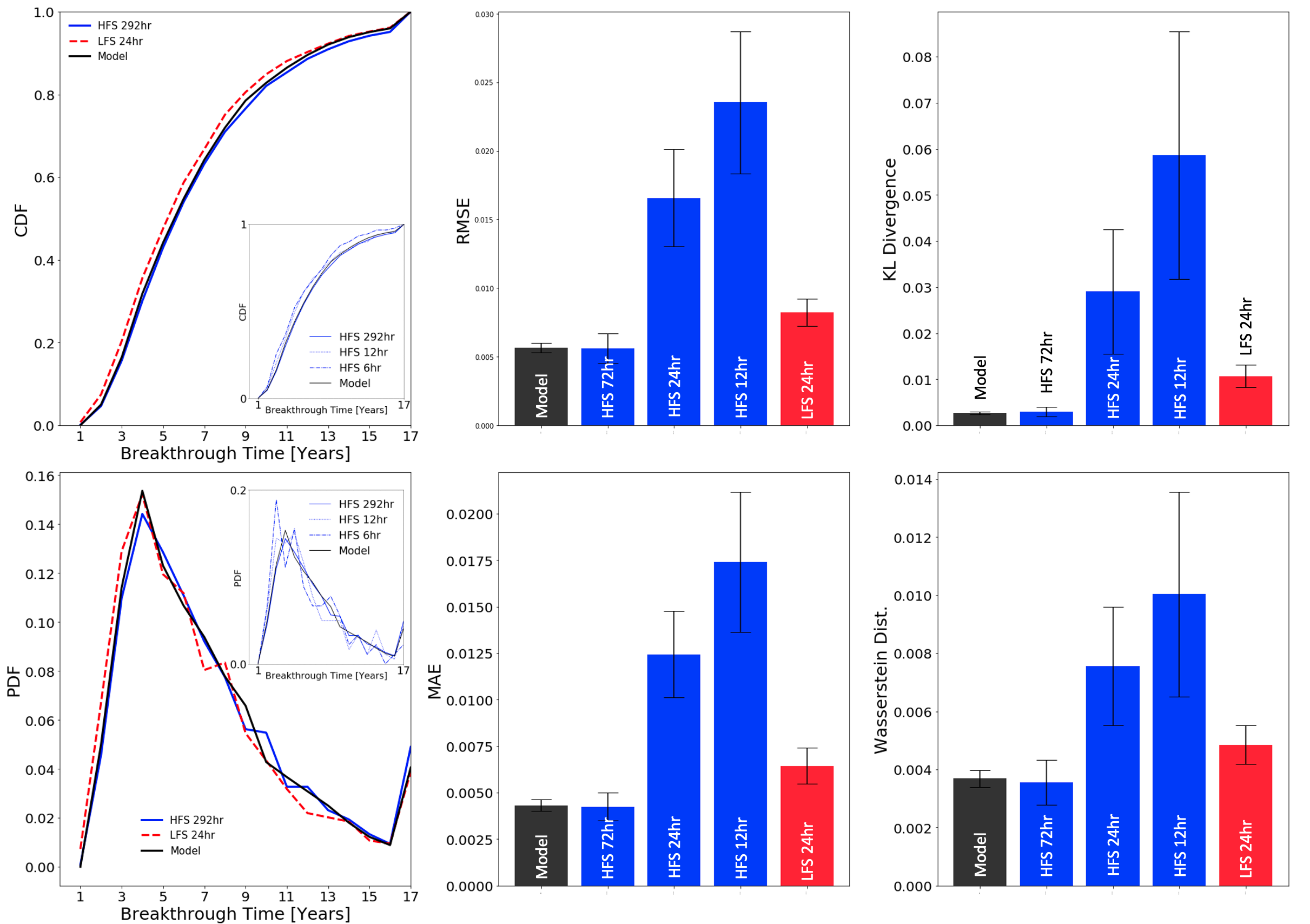}
\end{center}
\caption{Left: The converged CDF (top) and PDF (bottom) of breakthrough time is calculated using MC simulations of HFS, LFS, and the CNN surrogate model. The CDF and PDF calculated from varying amounts of HFS are displayed on the subplots. Bar plots: RMSE (Middle-top), MAE (Middle-bottom), KL Divergence (Right-top), and Wasserstein distance (Right-bottom) from PDF calculated using CNN model, HFS, and LFS.}
\label{fig: results UQ}
\end{figure}

Figure~\ref{fig: results UQ} exhibits the CDF and PDF of $T_\text{break}$ alternatively computed with HFS and LFS Monte Carlo and with the CNN trained on the multi-fidelity data. The distributions obtained via Monte Carlo consisting of 282 hours of HFS are treated as ground truth. The distributions obtained from 24 hours of LFS involve a sufficient number of samples for the error to be attributable solely to the low resolution, i.e., to the disretization errors in solving PDEs.
%
%The CDF/PDF from 282 hours of HFS and 24 hours of LFS represent the converged MC simulation solutions. 
%We take the 282 hours of HFS MC simulations to be the ground truth. 
%The discretization error of the LFS is evident in Fig~\ref{fig: results UQ}, as the converged solutions of LFS do not match that of HFS. 
The numbers of HFS samples generated during either 6 or 12 hours of simulations are insufficient for Monte Carlo to converge, leading to the appreciable errors in estimation of PDF and CDF of $T_\text{break}$. The CNN trained on multi-fidelity  data yields accurate estimates of these quantities, while requiring only 12 hours of data generation. 

In addition to visual comparison, the alternative strategies for estimation of the distributions of $T_\text{break}$ are compared in terms of RMSE, mean absolute error (MAE), the Kullback-Leibler (KL) divergence, and the Wasserstein distance. The UQ task was repeated 50 times, with Fig.~\ref{fig: results UQ} displaying the mean and standard deviation of these measures of discrepancy. We found 3200 forward passes of the CNN to be sufficient for the CDF/PDF estimates to converge; this UQ task took about 10 minutes, whereas an equivalent HFS Monte Carlo takes 194 hours. By every discrepancy measure, the CNN estimates outperform the converged LFS Monte Carlo and are at least as accurate as the HFS Monte Carlo using 72 hours of data. Likewise, the CNN estimates are vastly more accurate than the HFS Monte Carlo of a similar data-generation budget.

%%%%%%%%%%
\section{Conclusions}
\label{sec: Conclusions}
We proposed a transfer learning-based approach to train a CNN on multi-fidelity (e.g., multi-resolution) data. High- and low-fidelity images were generated by solving a PDE on fine and coarse meshes, respectively. The performance of our algorithm was tested on a system of nonlinear parabolic PDEs governing multi-phase flow in a heterogeneous porous medium with uncertain (random) permeability. A quantity of interest (QoI) in this example is PDF or CDF of the breakthrough time of an invading fluid. %This method was tested with various ratios of HFS and LFS for a given data-generation budget. The trained CNN was made to calculate distributions of a quantity of interest, and the distributions were compared to those generated using MC simulations. %The CNN performance was assessed in terms of RMSE on test data. The trained networks were used to model multi-phase fluid displacement and calculate the PDF and CDF of break-through times based on random permeability fields. The PDF and CDF of the models were compared to that of MC simulations using RMSE, MAE, KL divergence and Wasserstein distance. 
Our analysis leads to the following major conclusions.

\begin{enumerate}
    \item CNN surrogates trained on multi-fidelity data provide an accurate approximation of the PDE solution on the fine mesh, even for highly nonlinear problems that exhibit sharp dynamic fronts. A forward pass of the CNN surrogate is two orders of magnitude faster than a PDE solution on the fine-mesh. This speed-up makes CNN surrogates an invaluable tool for ensemble-based computation of the PDF/CDF of a QoI.
    
    \item CNN training on multi-fidelity data reduces the data-generation budget 7-fold relative to to CNN training on HFS data alone. If the budget is relatively small, the CNN trained on the LFS data is more accurate than its HFS-trained counterpart.  As the budget increases, the opposite is true. This finding is reminiscent of the cost-constrained selection between high- and low-fidelity models in the context of ensemble-based simulations. 
    
    \item For a small data-generation budget (12 hours, in our example), the CNN trained on multi-fidelity data exhibits an appreciably smaller RMSE on test data than the CNNs trained on either HFS or LFS data. Performance of the multi-fidelity CNN depends on the ratio between HFS and LFS in the training set. Theoretical studies from multilevel Monte Carlo can be used to guide the selection of an optimal mix of low- and high-fidelity data.
    
%    \item The CNN performance depends of ratio between the HFS and LFS data is critical when using multiple fidelity data to train a CNN.  In our problem, having 573 LFS and 100 HFS (12 hours data-generation budget) yielded results similar to CNNs trained on 79 hours of HFS data. However, training the CNN using 1 LFS and 197 HFS (12 hours data-generation budget) yielded worse results than just using 12 hours of HFS data.
    
    \item The CNN trained on multi-fidelity data is largely insensitive to the discretization error of LFS. CNN-derived estimates of the PDF and CDF of the QoI are close to those of converged high-fidelity Monte Carlo; but the former are three orders of magnitude faster to obtain than the latter. %there is a measurable difference between the distributions generated using the trained CNN and converged low-scale MC simulations.
    
%    \item CNN trained on multiple-fidelity is applicable for UQ tasks. The UQ performance of the CNN trained using 12 total hours of multiple-fidelity data was similar to that of 72 hours of high-fidelity MC simulations; the CNN based recreation of the PDF/CDF was more accurate than those generated using 12 or 24 hours of high-fidelity MC simulations.
\end{enumerate}

%% The Acknowledgements part is started with the command \acknowledgements;
%% acknowledgements are then done as normal sections before appendix
%% \acknowledgements

\acknowledgements

This research was supported in part by Air Force Office of Scientific Research under award number FA9550-18-1-0474; by the Advanced Research Projects Agency-Energy (ARPA-E), U.S. Department of Energy, under Award Number  DE-AR0001202; and by a gift from Total.

\begin{figure}[!h]
\begin{center}
\includegraphics[width=1\textwidth]{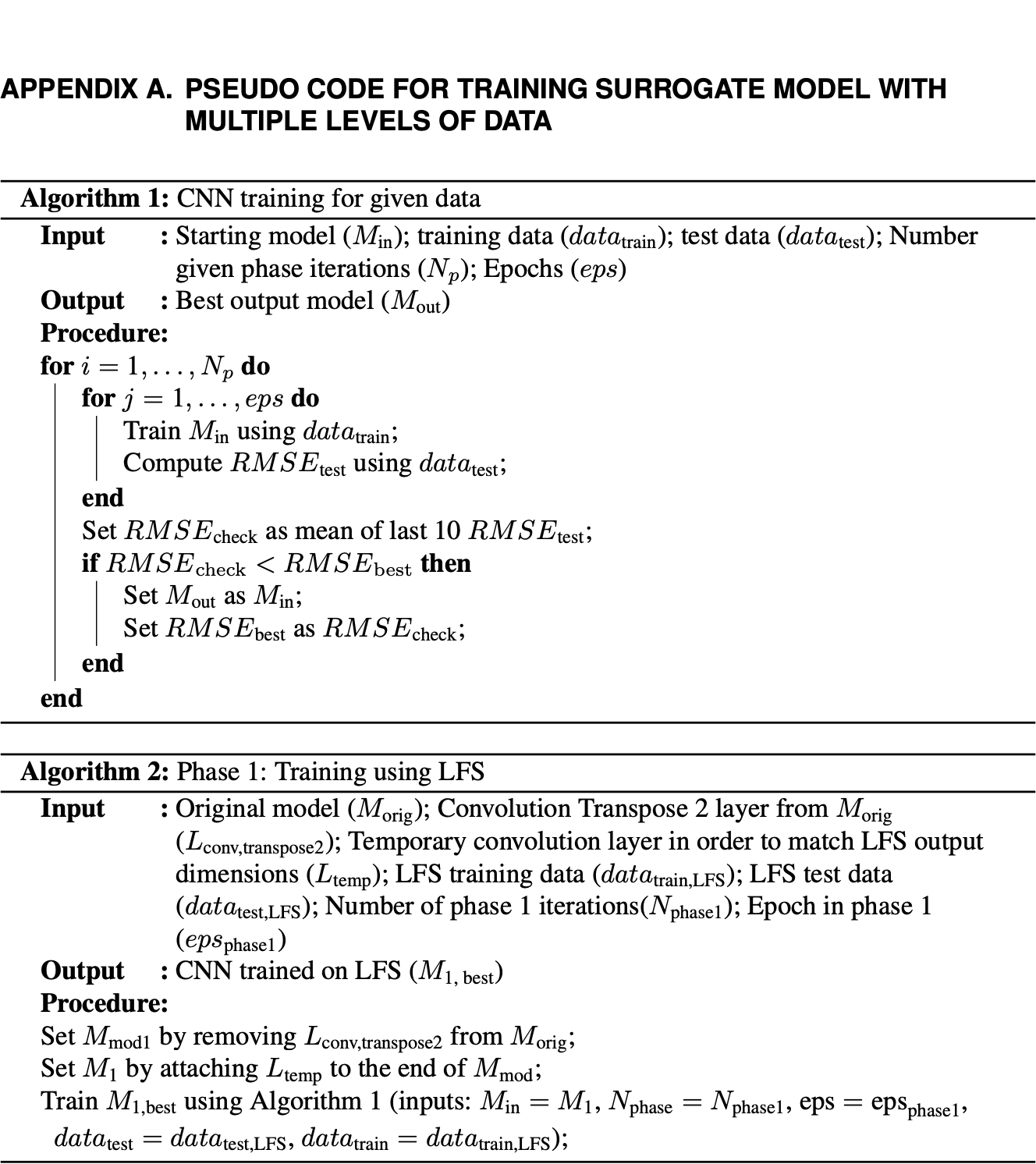}
\end{center}
\end{figure}

\begin{figure}[!h]
\begin{center}
\includegraphics[width=1\textwidth]{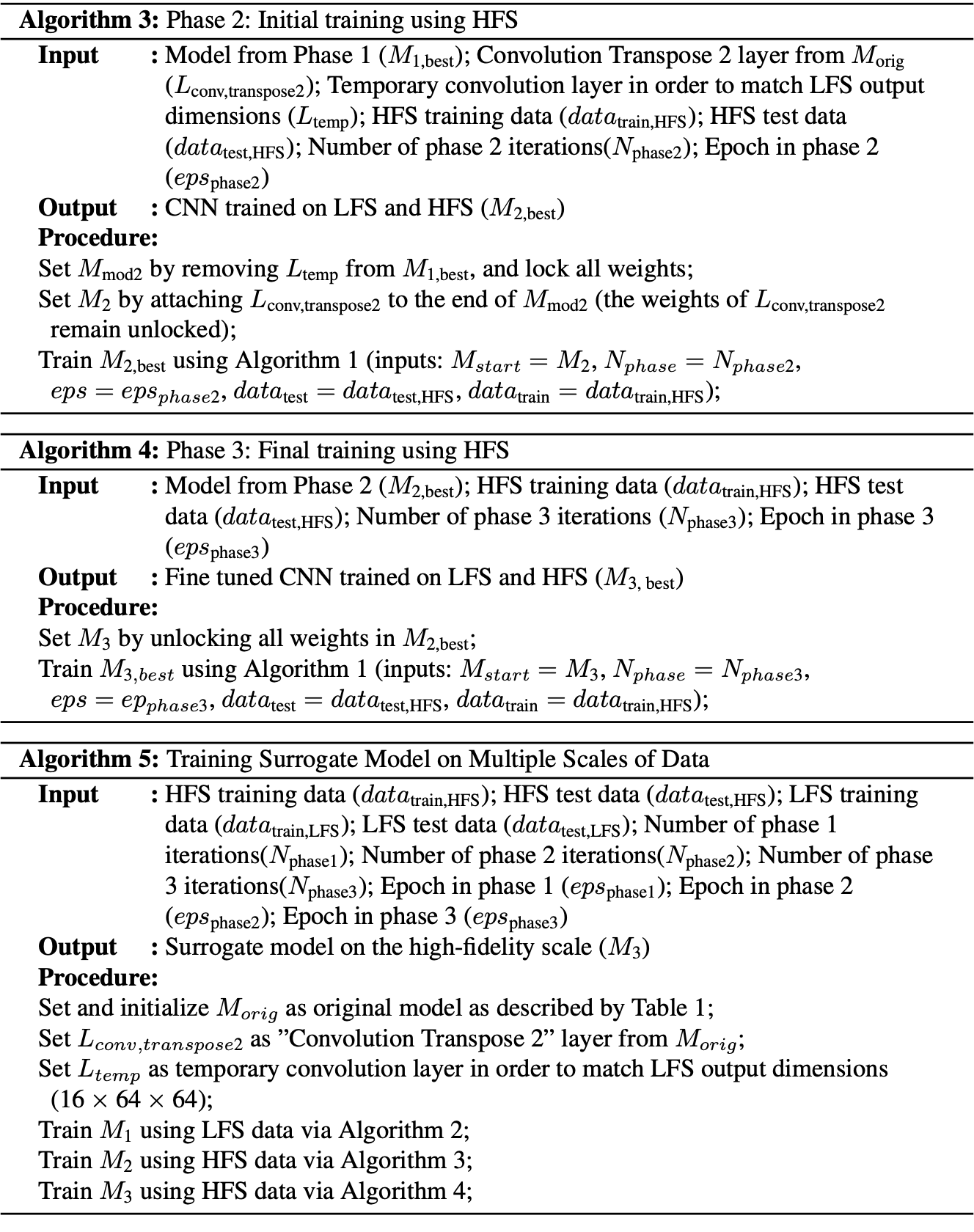}
\end{center}
\end{figure}

\pagebreak

\bibliographystyle{Bibliography_Style}

\bibliography{multifidelityCNN}

\begin{thebibliography}{37}
\expandafter\ifx\csname natexlab\endcsname\relax\def\natexlab#1{#1}\fi
\expandafter\ifx\csname url\endcsname\relax
  \def\url#1{\texttt{#1}}\fi
\expandafter\ifx\csname urlprefix\endcsname\relax\def\urlprefix{URL }\fi

\bibitem[{Appleyard et~al.(1981)Appleyard, Cheshire, and
  Pollard}]{appleyard1981special}
Appleyard, J.R., Cheshire, I.M., and Pollard, R.K., \titlecap{Special
  techniques for fully implicit simulators}, {\em Proceedings of the European
  Symposium on Enhanced Oil Recovery}, Bournemouth, UK, pp.~395--408, 1981.

\bibitem[{Aziz(1979)}]{aziz1979petroleum}
Aziz, K., {\em \titlecap{Petroleum reservoir simulation}}, Vol. 476, Applied
  Science Publishers, New York, 1979.

\bibitem[{Boso and Tartakovsky(2018)}]{boso-2018-information}
Boso, F. and Tartakovsky, D.M., \titlecap{Information-theoretic approach to
  bidirectional scaling}, {\em Water Resour. Res.}, vol.~{\bf 54}, no.~7,
  pp.~4916--4928, 2018.

\bibitem[{Breiman(2001)}]{Breiman-2014-Random}
Breiman, L., \titlecap{Random forests}, {\em Mach. Learn.}, vol.~{\bf 45},
  no.~1, pp.~5--32, 2001.

\bibitem[{Corey(1954)}]{corey1954interrelation}
Corey, A.T., \titlecap{The interrelation between gas and oil relative
  permeabilities}, {\em Producers Month.}, vol.~{\bf 19}, no.~1, pp.~38--41,
  1954.

\bibitem[{Couckuyt et~al.(2014)Couckuyt, Dhaene, and
  Demeester}]{Couckuyt-2014-ooDACE}
Couckuyt, I., Dhaene, T., and Demeester, P., \titlecap{{SooDACE} Toolbox: {A}
  Flexible Object-Oriented Kriging Implementation}, {\em J. Mach. Learn. Res.},
  vol.~{\bf 15}, pp.~3183--3186, 2014.

\bibitem[{Donahue et~al.(2014)Donahue, Jia, Vinyals, Hoffman, Zhang, Tzeng, and
  Darrell}]{donahue2014decaf}
Donahue, J., Jia, Y., Vinyals, O., Hoffman, J., Zhang, N., Tzeng, E., and
  Darrell, T., \titlecap{Decaf: A deep convolutional activation feature for
  generic visual recognition}, {\em International Conference on Machine
  Learning}, PMLR, pp.~647--655, 2014.

\bibitem[{Durlofsky(2005)}]{durlofsky2005upscaling}
Durlofsky, L.J., \titlecap{Upscaling and gridding of fine scale geological
  models for flow simulation}, {\em 8th International Forum on Reservoir
  Simulation}, Vol. 2024, Iles Borromees, Stresa, Italy, pp.~1--59, 2005.

\bibitem[{Friedman et~al.(2001)Friedman, Hastie, and
  Tibshirani}]{friedman2001elements}
Friedman, J., Hastie, T., and Tibshirani, R., {\em \titlecap{The elements of
  statistical learning}}, Vol.~1, Springer, New York, 2001.

\bibitem[{Fuks and Tchelepi(2020)}]{fuks2020}
Fuks, O. and Tchelepi, H., \titlecap{Limitations of physics informed machine
  learning for nonlinear two-phase transport in porous media}, {\em J. Mach.
  Learn. Model. Comput.}, vol.~{\bf 1}, no.~1, pp.~19--37, 2020.

\bibitem[{Giles(2008)}]{giles2008multilevel}
Giles, M.B., \titlecap{Multilevel Monte Carlo path simulation}, {\em Oper.
  Res.}, vol.~{\bf 56}, no.~3, pp.~607--617, 2008.

\bibitem[{Haghighat et~al.(2021)Haghighat, Raissi, Moure, Gomez, and
  Juanes}]{haghighat2021physics}
Haghighat, E., Raissi, M., Moure, A., Gomez, H., and Juanes, R., \titlecap{A
  physics-informed deep learning framework for inversion and surrogate modeling
  in solid mechanics}, {\em Comput. Meth. Appl. Mech. Engrg.}, vol.~{\bf 379},
  p.~113741, 2021.

\bibitem[{Heinrich(1998)}]{heinrich1998monte}
Heinrich, S., \titlecap{Monte Carlo complexity of global solution of integral
  equations}, {\em J. Complexity}, vol.~{\bf 14}, no.~2, pp.~151--175, 1998.

\bibitem[{Heinrich(2001)}]{heinrich2001multilevel}
Heinrich, S., \titlecap{Multilevel Monte Carlo methods}, {\em International
  Conference on Large-Scale Scientific Computing}, Springer, pp.~58--67, 2001.

\bibitem[{Hwang and Martins(2018)}]{Hwang-2018-fast}
Hwang, J.T. and Martins, J.R.R.A., \titlecap{A fast-prediction surrogate model
  for large datasets}, {\em Aerospace Sci. Tech.}, vol.~{\bf 75}, pp.~74--87,
  2018.

\bibitem[{Jiang and Learned-Miller(2017)}]{jiang2017face}
Jiang, H. and Learned-Miller, E., \titlecap{Face detection with the faster
  R-CNN}, {\em 2017 12th IEEE International Conference on Automatic Face \&
  Gesture Recognition (FG 2017)}, IEEE, pp.~650--657, 2017.

\bibitem[{Karpathy and Fei-Fei(2015)}]{karpathy2015deep}
Karpathy, A. and Fei-Fei, L., \titlecap{Deep visual-semantic alignments for
  generating image descriptions}, {\em Proceedings of the IEEE Conference on
  Computer Vision and Pattern Recognition}, pp.~3128--3137, 2015.

\bibitem[{Kingma and Ba(2014)}]{kingma2014adam}
Kingma, D.P. and Ba, J., \titlecap{Adam: A method for stochastic optimization},
  {\em arXiv:1412.6980}, 2014.

\bibitem[{Lagaris et~al.(1998)Lagaris, Likas, and
  Fotiadis}]{lagaris1998artificial}
Lagaris, I.E., Likas, A., and Fotiadis, D.I., \titlecap{Artificial neural
  networks for solving ordinary and partial differential equations}, {\em IEEE
  Trans. Neural Networks}, vol.~{\bf 9}, no.~5, pp.~987--1000, 1998.

\bibitem[{Lee and Kang(1990)}]{lee1990neural}
Lee, H. and Kang, I.S., \titlecap{Neural algorithm for solving differential
  equations}, {\em J. Comput. Phys.}, vol.~{\bf 91}, no.~1, pp.~110--131, 1990.

\bibitem[{Mo et~al.(2019{\natexlab{a}})Mo, Zabaras, Shi, and Wu}]{mo2019deep_b}
Mo, S., Zabaras, N., Shi, X., and Wu, J., \titlecap{Deep autoregressive neural
  networks for high-dimensional inverse problems in groundwater contaminant
  source identification}, {\em Water Resour. Res.}, vol.~{\bf 55}, no.~5,
  pp.~3856--3881, 2019{\natexlab{a}}.

\bibitem[{Mo et~al.(2019{\natexlab{b}})Mo, Zhu, Zabaras, Shi, and
  Wu}]{mo2019deep_a}
Mo, S., Zhu, Y., Zabaras, N., Shi, X., and Wu, J., \titlecap{Deep convolutional
  encoder-decoder networks for uncertainty quantification of dynamic multiphase
  flow in heterogeneous media}, {\em Water Resour. Res.}, vol.~{\bf 55}, no.~1,
  pp.~703--728, 2019{\natexlab{b}}.

\bibitem[{Montgomery and Evans(2018)}]{Montgomery-2018-Second}
Montgomery, D.C. and Evans, D.M., \titlecap{Second-order response surface
  designs in computer simulation}, {\em Aerospace Sci. Tech.}, vol.~{\bf 75},
  pp.~74--87, 2018.

\bibitem[{M\"uller et~al.(2013)M\"uller, Jenny, and
  Meyer}]{muller-2013-Multifidelity}
M\"uller, F., Jenny, P., and Meyer, D.W., \titlecap{Multilevel {Monte~Carlo}
  for two phase flow and {Buckley-Leverett} transport in random heterogeneous
  porous media}, {\em J. Comput. Phys.}, vol.~{\bf 250}, pp.~685--702, 2013.

\bibitem[{Paleologos et~al.(1996)Paleologos, Neuman, and
  Tartakovsky}]{paleologos-1996-effective}
Paleologos, E.K., Neuman, S., and Tartakovsky, D.M., \titlecap{Effective
  hydraulic conductivity of bounded, strongly heterogeneous porous media}, {\em
  Water Resour. Res.}, vol.~{\bf 32}, no.~5, pp.~1333--1341, 1996.

\bibitem[{Paszke et~al.(2019)Paszke, Gross, Massa, Lerer, Bradbury, Chanan,
  Killeen, Lin, Gimelshein, Antiga, Desmaison, Kopf, Yang, DeVito, Raison,
  Tejani, Chilamkurthy, Steiner, Fang, Bai, and Chintala}]{NEURIPS2019_9015}
Paszke, A., Gross, S., Massa, F., Lerer, A., Bradbury, J., Chanan, G., Killeen,
  T., Lin, Z., Gimelshein, N., Antiga, L., Desmaison, A., Kopf, A., Yang, E.,
  DeVito, Z., Raison, M., Tejani, A., Chilamkurthy, S., Steiner, B., Fang, L.,
  Bai, J., and Chintala, S., 2019. Pytorch: An imperative style,
  high-performance deep learning library, {\em Advances in Neural Information
  Processing Systems 32}. H.~Wallach, H.~Larochelle, A.~Beygelzimer,
  F.~d\textquotesingle Alch\'{e}-Buc, E.~Fox, and R.~Garnett, Eds. Curran
  Associates, Inc., pp.~8024--8035.

\bibitem[{Peherstorfer(2019)}]{Peherstorfer-2019-Multifidelity}
Peherstorfer, B., \titlecap{Multifidelity {Monte Carlo} Estimation with
  Adaptive Low-Fidelity Models}, {\em SIAM/ASA J. Uncert. Quant.}, vol.~{\bf
  7}, no.~2, p.~579–603, 2019.

\bibitem[{Raissi et~al.(2019)Raissi, Perdikaris, and
  Karniadakis}]{raissi2019physics}
Raissi, M., Perdikaris, P., and Karniadakis, G.E., \titlecap{Physics-informed
  neural networks: A deep learning framework for solving forward and inverse
  problems involving nonlinear partial differential equations}, {\em J. Comput.
  Phys.}, vol.~{\bf 378}, pp.~686--707, 2019.

\bibitem[{Sinsbeck and Tartakovsky(2015)}]{sinsbeck-2015-impact}
Sinsbeck, M. and Tartakovsky, D.M., \titlecap{Impact of data assimilation on
  cost-accuracy tradeoff in multifidelity models}, {\em SIAM/ASA J. Uncert.
  Quant.}, vol.~{\bf 3}, no.~1, pp.~954--968, 2015.

\bibitem[{Tang et~al.(2020)Tang, Liu, and Durlofsky}]{tang2020deep}
Tang, M., Liu, Y., and Durlofsky, L.J., \titlecap{A deep-learning-based
  surrogate model for data assimilation in dynamic subsurface flow problems},
  {\em J. Comput. Phys.}, p.~109456, 2020.

\bibitem[{Tartakovsky and Neuman(1998)}]{tartakovsky-1998-transient-2}
Tartakovsky, D.M. and Neuman, S.P., \titlecap{Transient effective hydraulic
  conductivities under slowly and rapidly varying mean gradients in bounded
  three-dimensional random media}, {\em Water Resour. Res.}, vol.~{\bf 34},
  no.~1, pp.~21--32, 1998.

\bibitem[{Taverniers et~al.(2020)Taverniers, Bosma, and
  Tartakovsky}]{taverniers2020accelerated}
Taverniers, S., Bosma, S.B., and Tartakovsky, D.M., \titlecap{Accelerated
  multilevel Monte Carlo with kernel-based smoothing and Latinized
  stratification}, {\em Water Resour. Res.}, vol.~{\bf 56}, no.~9,
  p.~e2019WR026984, 2020.

\bibitem[{Tripathy and Bilionis(2018)}]{tripathy2018deep}
Tripathy, R.K. and Bilionis, I., \titlecap{Deep UQ: Learning deep neural
  network surrogate models for high dimensional uncertainty quantification},
  {\em J. Comput. Phys.}, vol.~{\bf 375}, pp.~565--588, 2018.

\bibitem[{Xiu(2010)}]{xiu-2010-numerical}
Xiu, D., {\em \titlecap{Numerical Methods for Stochastic Computations: A
  spectral Methd Approach}}, Princeton University Press, Princeton, NJ, 2010.

\bibitem[{Yang et~al.(2020)Yang, Wang, and Tartakovsky}]{yang-2020-resource}
Yang, L., Wang, P., and Tartakovsky, D.M., \titlecap{Resource-constrained model
  selection for uncertainty propagation and data assimilation}, {\em SIAM/ASA
  J. Uncert. Quant.}, vol.~{\bf 8}, no.~3, pp.~1118--1138, 2020.

\bibitem[{Zhou and Tartakovsky(2021)}]{zitong2020}
Zhou, Z. and Tartakovsky, D.M., \titlecap{Markov chain {Monte~Carlo} with
  neural network surrogates: {Application} to contaminant source
  identification}, {\em Stoch. Environ. Res. Risk Assess.}, vol.~{\bf 35},
  no.~3, pp.~639--651, 2021.

\bibitem[{Zhu et~al.(2019)Zhu, Zabaras, Koutsourelakis, and
  Perdikaris}]{zhu2019physics}
Zhu, Y., Zabaras, N., Koutsourelakis, P.S., and Perdikaris, P.,
  \titlecap{Physics-constrained deep learning for high-dimensional surrogate
  modeling and uncertainty quantification without labeled data}, {\em J.
  Comput. Phys.}, vol.~{\bf 394}, pp.~56--81, 2019.

\end{thebibliography}
\end{document}